\documentclass[fleqn,leqno]{article}
\usepackage{amssymb,latexsym}
\usepackage{amsfonts}
\usepackage[noload]{qtree}
 \usepackage[dvips]{color}

\usepackage[pdftex=true]{hyperref}
\hypersetup{
 bookmarks=true, colorlinks=flase, breaklinks=true, pdfstartview={FitH}, unicode=true }
\newcommand{\bi}{\bigskip}
\newcommand{\sm}{\smallskip}

\newcommand{\wh}{\widehat}
\newcommand{\ee}{\end{equation}}
\newcommand{\eea}{\end{eqnarray}}
\newcommand{\bean}{\begin{eqnarray*}}
\newcommand{\eean}{\end{eqnarray*}}
\newif\ifpctex
\newcommand{\noi}{\noindent}

\newcommand{\ve}{\varepsilon}
\newcommand{\wt}{\widetilde}

\newtheorem{theorem}{Theorem}
\newtheorem{remark}{Remark}
\newtheorem{example}{Example}
\newcommand{\D}{\displaystyle}

\newcommand{\qad}{$\qquad \square$}
\newtheorem{proposition}{Proposition}[section]
\newtheorem{definition}[proposition]{Definition}
\newtheorem{lemma}[proposition]{Lemma}
\newtheorem{corollary}[proposition]{Corollary}

\newtheorem{conjecture}[proposition]{Conjecture}
\renewcommand{\theequation}{\thesection.\arabic{equation}}

 \newcommand{\rand}[1]{}

\newcommand{\La}{\Longrightarrow}

\newcommand{\Rand}[1]{}

\newcommand{\be}[1]{\Rand{\vspace{0,6cm}\tt #1}\begin{equation}\label{#1}}
\newcommand{\bea}[1]{\Rand{\vspace{0,7cm}\tt
#1\vspace{-0,7cm}}\begin{eqnarray}\label{#1}} \marginparwidth2.5cm

\newcommand{\beL}[1]{\Rand{\vspace{0,6cm}\tt #1}\begin{lemma}\label{#1}}
\newcommand{\beD}[1]{\Rand{\vspace{0,6cm}\tt #1}\begin{definition}\label{#1}}
\newcommand{\beT}[1]{\Rand{\vspace{0,6cm}\tt #1}\begin{theorem}\label{#1}}
\newcommand{\beP}[1]{\Rand{\vspace{0,6cm}\tt #1}\begin{proposition}\label{#1}}
\newcommand{\beC}[1]{\Rand{\vspace{0,6cm}\tt #1}\begin{corollary}\label{#1}}
\newcommand{\beCj}[1]{\Rand{\vspace{0,6cm}\tt #1}\begin{conjecture}\label{#1}}

\newcommand{\ep}{\end{proposition}}

\newcommand{\suml}{\sum\limits}
\newcommand{\intl}{\int\limits}
\newcommand{\liml}{\lim\limits}

\newcommand{\limsupl}{\limsup\limits}

\newcommand{\R}{\mathbb{R}}
\newcommand{\N}{\mathbb{N}}

\newcommand{\I}{\mathbb{I}}

\newcommand{\Q}{\mathbb{Q}}
\newcommand{\text}{\mbox}

\newcommand{\CE}{{\mathcal E}}
\newcommand{\CF}{{\mathcal F}}

\newcommand{\CL}{{\mathcal L}}

\newcommand{\CM}{{\mathcal M}}

\newcommand{\CP}{{\mathcal P}}
\newcommand{\CU}{{\mathcal U}}
 
\newcommand{\CW}{{\mathcal W}}

\newcommand{\tto}{{_{\D \Longrightarrow \atop t \to \infty}}}

\newcommand{\ttooo}{{_{\D \longrightarrow \atop t \to -\infty}}}

\newcommand{\Nto}{{_{\D \Longrightarrow \atop N \to \infty}}}

\newcommand{\Ntoo}{{_{\D \longrightarrow \atop N \to \infty}}}

\newcommand{\la}{\longrightarrow}

\begin{document}

% Wenn rand-Notizen wegfallen sollen: bei \renewcommand\Rand die %%%
%l"oschen!

\title{Multiscale analysis: Fisher-Wright diffusions with
rare mutations and selection, logistic branching system}
\author{{\bf Donald A. Dawson$^{1,2}$ \mbox{ } Andreas Greven$^{2,3}$}}

\date{\small  \today }%\\ {\tt script/dg/gaertner15.tex}}
\maketitle

\begin{abstract}
We study two types of stochastic processes, a mean-field spatial system of interacting Fisher-Wright diffusions with
an inferior and an advantageous type with rare mutation (inferior to advantageous) and
a (mean-field) spatial system of supercritical branching random walks with an additional
deathrate which is quadratic in the local number of particles.
The former describes a standard two-type population under selection, mutation,
the latter models describe a population under scarce resources causing additional
death at high local population intensity. Geographic space is modelled
by $\{1, \cdots, N\}$.
The first process starts in an initial state with only the inferior
type present or an exchangeable configuration and the second one with
a single initial particle.
{This material is a special case of the theory developed in \cite{DGsel}.}

We study the behaviour in two time windows, first between time 0 and $T$ and
secondly after a large time when in the Fisher-Wright model
the rare mutants succeed respectively in the branching random walk
the particle population reaches a positive
spatial intensity. It is shown that the second phase for both models sets in after
time $\alpha^{-1} \log N$, if $N$ is the size of geographic space and $N^{-1}$
the rare mutation rate and $\alpha \in (0, \infty)$  depends on the other parameters.
We identify the limit
dynamics in both time windows and for both models as a nonlinear Markov dynamic
(McKean-Vlasov dynamic) respectively a corresponding random entrance law from
time $-\infty$ of this dynamic.

Finally we explain that the two processes are just two sides of the very same
coin, a fact arising from duality, in particular the particle model  generates
the genealogy of the Fisher-Wright diffusions with selection and mutation.
We discuss the extension of this duality relation to a multitype model with
more than two types.

\end{abstract}


\vspace{2cm}

\footnoterule
\noi
\hspace*{0.3cm}{\footnotesize $^1$ {School of Mathematics and Statistics,
Carleton University, Ottawa K1S 5B6, Canada, e-mail:
ddawson@math.carleton.ca}}\\
\hspace*{0.3cm}{\footnotesize $^2$
{Research supported by NSERC, DFG-NWO
Forschergruppe 498.}}\\
\hspace*{0.3cm}{\footnotesize $^3$
{Department Mathematik, Universit\"at Erlangen-N\"urnberg,
Bismarckstra{\ss}e 1 1/2,
D-91054 Erlangen, Germany, e-mail: greven@mi.uni-erlangen.de}}
\newpage

    \tableofcontents

 \newpage

		    \newcounter{secnum}
		    \setcounter{secnum}{\value{section}}
		    \setcounter{section}{-1}
		    %\addtocontents{toc}{Introduction}
		    \setcounter{equation}{0}
		    \renewcommand{\theequation}{\mbox{\arabic{secnum}.\arabic{equation}}}
		    \section{Introduction}
\label{s.intro}

\subsection{Motivation and background}
\label{ss.motive}

We study here features of the longtime behavior of two models for the stochastic evolution of populations,
the classical (mean-field) spatial version of a system of interacting Fisher-Wright
diffusions with selection and mutation on the one hand and a logistic spatial
branching particle model both on geographic space $\{1, \cdots, N\}$.
We shall later explain the mathematical relation between them.
\sm

{\em Fisher-Wright model}. This process comes from population genetics and
 models a population of individuals of two types evolving under
migration, resampling, selection and mutation.
It is the many individual limit of a discrete model.
Migration here means individuals move in geographic space, resampling that
pairs are replaced by an offspring pair each choosing a parent at random
and adopting the parents type, mutation is a spontaneous change of type
of an individual
and under selection the choice of the parent in the resampling event
is biased according to the parents' fitness.

Here we are particularly
concerned with a situation where we have an inferior type and an advantageous
type. The case we are interested in is that the mutation rate from
inferior to advantageous is very small so
that in finite time we expect $O(1)$-many mutations in all of space as space gets large.
We want to follow the population through the {\em emergence} and {\em fixation} of the
whole population in the advantageous type.
\sm

{\em Logistic branching model.} Here we consider a population of particles
which migrate in space, have offspring at a certain rate $s$, die with a certain
rate  (all particles act independently) but here the risk of
death is at a rate increasing with the population size of the site and being zero
for one particle. The latter
mechanism induces an interaction of the families, which therefore do not
evolve anymore independently of each other.

Here for this model we want to see how the population {starting from one particle} spreads and
eventually colonizes the whole space as time evolves, meaning that a positive spatial intensity
is reached and a local equilibrium situation arises where locally the process neither
becomes extinct nor grows and becomes infinitely large as $t \to \infty$.
This is in contrast to the behaviour of classical branching models with their
survival versus extinction dichotomy in finite geographic space and
reflects the limited resources in a given colony.
\sm

{\em McKean-Vlasov equation and random entrance laws.}
The techniques we use to study the questions raised above for two models is the mean-field
limit, where we choose migration to occur according to the uniform distribution
on $\{1,\cdots, N\}$ and rare mutation having rate $N^{-1}$
and where we let $N \to\infty$. As limit dynamics a nonlinear Markov
process of the type first introduced by McKean \cite{Mc} arises, i.e. an evolution with a generator where the parameters depend on the current
state of the process {\em and} on the current law. The transition probabilities
solve the {\em McKean-Vlasov equation} which has a similar structure to the equations introduced by Vlasov to describe the dynamics of a plasma consisting of charged particles with long-range interaction.
A similar scenario holds for the branching particle system.

The first rigorous and systematic analysis of mean-field models and resulting
McKean-Vlasov dynamics is  G\"artner's fundamental paper in 1988
\cite{Gar} which established the existence and uniqueness of weak solutions for a general class of McKean-Vlasov equations and the associated non-linear martingale problems but under the condition that the diffusion matrix is {\em
strictly} positive definite.  Hence it does not cover the  case dealt with in this paper due to the fact that for the Fisher Wright diffusions the diffusion function vanishes at  the boundary.

% G\"{a}rtner \cite{Gar}
%{\bf xxx1 History, J\"urgen's work}

We must also extend here the methodology of the McKean-Vlasov equation
further in order to describe the
limiting behaviour in different time windows.  In order to also consider a late time window we need first of all to introduce the notion of an entrance law from
time $-\infty$ (our time parameter has the form $T_N+t, \quad t\in \R$
and $T_N \to \infty$ as $N\to \infty$) but since in the initial time phase
some randomness is involved (rare mutation respectively very small early particle
intensity) we even have to work with {\em random entrance laws} to the
McKean-Vlasov equation.

{\em Duality}.
The mathematical structure which relates our two models is the fact that
they are in duality, i.e. expectations of certain functionals under
one dynamic are given by expectations of appropriate functionals under
the other dynamic with the time direction reversed. This will be explained in detail later on.

The duality relation which we present here is a special case of a broader new
duality theory for multitype models, which allows a historical and
genealogical interpretation and which has been developed in
\cite{DGsel} covering much more general situations than we
can discuss here. {More information on the particle system can be
found in \cite{Schirm10}.}

\bi

\begin{remark}
In the framework just sketched we can analyse the features of the population as
described above. If one wants to adapt a more realistic model for
geographic space the mean-field limit has to be replaced by the
{\em hierarchical mean-field limit} for which the present analysis is a key
ingredient. Then it is possible to study
asymptotically two-dimensional geographic space via its approximations
by the hierarchical group of order $N$ and $N \to \infty$. This also allows us to investigate
the question of {\em universality} of the behaviour. All this is carried out in
\cite{DGsel} and we refer the interested reader to this paper.
\end{remark}

\subsection{Outline}
\label{ss.outline}

In section \ref{s.FWraremut} we shall present the Fisher-Wright model, in section
\ref{s.logisticbrw} the logistic branching model and in section \ref{s.dualrel} the
connection between both via duality.

\section[The Fisher-Wright model with rare mutation and selection]
{The Fisher-Wright model with rare mutation and selection: Behaviour in two time windows}

\label{s.FWraremut}

Here we introduce the Fisher-Wright model, some relevant time windows,
state its properties in two time
windows and consider the McKean-Vlasov limit for the model in five separate subsections.

\subsection{A two-type mean-field diffusion model and its description}
\label{ss.2mfdiff}

We study a population with two types, one  of low and the other of high fitness, where
initially all the population is of the lower type but by a rare mutation
the advantageous type appears and spreads with time. The population is
described specifying the proportion of the type 1 (the inferior)  in every spatial
colony.

Formally we look at a process
$(X^N(t))_{t \geq 0}$ with $N \in \N$ of the following form:

\be{angr1}
X^N(t) =((x^N_1 (i,t),x^N_2 (i,t)); i=1,\dots,N),
\ee
\be{angr1b}
x^N_1 (i,0) =1, \quad i=1,\cdots,N,
\ee
satisfying the well-known SSDE:
\bea{angr2} d x^N_1 (i,t) = c(\bar x^N_1 (t) - x^N_1 (i,t)) dt
 &-& s\,x^N_1(i,t) x^N_2 (i,t) dt
  - \frac{m}{L} x^N_1 (i,t)dt \\
   &+& \sqrt{d \cdot x^N_1(i,t)x^N_2(i,t)} dw_1(i,t), \nonumber
    \eea

\bea{angr3}
d x^N_2 (i,t) = c(\bar x^N_2 (t) - x^N_2 (i,t)) dt
 &+& s\,x^N_1(i,t) x^N_2 (i,t) dt + \frac{m}{L} x^N_1 (i,t)dt \\
   &+& \sqrt{d \cdot x^N_1 (i,t)x^N_2 (i,t)} dw_2 (i,t),\nonumber
    \eea
    where $w_2 (i,t) = - w_1(i,t)$ and $\{(w_1(i,t))_{t \geq 0},
    i=1,\dots,N\}$ are i.i.d. Brownian motions, $m,d,s,L \in (0,\infty)$
and
\be{angr3b}
\bar {x}^N_\ell (t) = \frac{1}{N}\wh {x}^N_\ell (t)
\mbox{ with } \ell=1,2, \quad \wh x^N_\ell(t) = \suml^N_{i=1} x^N_\ell(i,t). \ee
Later we will use the parameter of mutation strength and size of
geographic space satisfying
\be{ad1}
L=N.
\ee

We can study this system in various ways, {\em locally} by looking at a
$K$ tagged sites
\be{angr5++}
(x^N_\ell(1,t), \cdots, x^N_\ell(K,t)), \quad \ell=1,2
\ee
or {\em globally} using the concept of the empirical measure of the
complete population:
\be{angr5} \Xi_N (t) := \frac{1}{N} \suml^N_{i=1}
\delta_{(x^N_1 (i,t), x^N_2(i,t)})\in
\mathcal{P}(\mathcal{P}(\{1,2\}))\ee and the empirical measure
process of either type:
\be{angr5+}
\Xi_N(t,\ell):=\frac{1}{N} \sum_{i=1}^N \delta_{x_\ell^N(i,t)}\in
\mathcal{P}([0,1]),\; \ell=1,2. \ee Note that since $x^N_1 (i,t) +
x^N_2 (i,t) =1$ it suffices in the case of two types to know one
component of the pair in (\ref{angr5+}), the other is then determined by this
condition. Then for two types we can  effectively replace $\CP(\CP(\{1,2\}))$ by $\CP([0,1])$
as states of the empirical measure.

However since very soon rare mutants appear somewhere in space they are present
for $t >0$, but of
course the set of sites where they are visible is {\em sparse}, we have to find a way to
describe this small subset where the advantageous type appears.
One (global) way is to consider the process
\be{angr3c}
\wh {x}^N_\ell (t) =  \suml^{N}_{i=1} x^N_\ell(i,t)
   \mbox{ with } \ell=1,2,
   \ee
but this way we lose the internal structure of the droplet of sites
with advantageous types of substantial mass.  In order to
keep track of the sparse set of sites at which nontrivial mass appears we
will give a random label to each site and define the following
{\em atomic-measure-valued process}.

 We assign independent of the process a point $a(j)$ randomly in $[0,1]$ to each site
 $j\in\{1,\dots,N\}$, that is, we define the collection
 \be{dd29a}
 \{a(j), \quad j=1,\cdots, N \} \mbox{  i.i.d. uniform on } [0,1].
 \ee
 We then associate with our process and a realization of the
 random labels a measure-valued process on $\CP( [0,1])$, which we denote by
 $(\gimel^{N,m}_t)_{t \geq 0}$ where
 \be{dd29}
 \gimel^{N,m}_t=\sum_{j=1}^N x^N_2(j,t) \delta_{a(j)}.
 \ee
The process $(\gimel^{N,m}_t)_{t \geq 0}$  describes the essential
features of the advantageous population.

\subsection{Two time windows for the spread of the advantageous type}
\label{ss.2tsadv}

The first question now is how the advantageous type develops in finite time
and the second on what time scale does the advantageous type take over
the whole population meaning we want to identify a time $T_{N,L}$ after
which the advantageous type has positive spatial intensity.
A key role is played by space and the fact that
the dynamic is random and not deterministic as in the scenario looked
at by \cite{Bu01}.
So we look at the cases $N=1$, $N$ large and $L$ small, $d=0,d>0$.
\sm

{\em Case $N=1$ (nonspatial)}

Here we find the time scales in which we emerge $(T_{N,L}=T_L)$ as $L \to \infty$ to be
\be{ad2}
O(s^{-1} \log L) \mbox{ for } d=0,
\ee
\be{ad3}
O(L)  \mbox{ for } d>0.
\ee
This qualitative difference is essentially  due to the fact that in the stochastic model the diffusion
can hit 0 by sheer random effects, while in the deterministic model the advantageous type
expands exponentially fast leading to a $\log L$ time scale.
\sm

{\em Case $N\gg 1$ (large spatial model)}

Here we take $L=N (T_{N,L}=T_N)$ and then we find also in the deterministic case the time scale
$s^{-1} \log N$, but also now in the stochastic model, $d>0$,  we have a time period
$\alpha^{-1} \log N$ but now the $\alpha$ will turn out to be strictly
smaller than $s$. In the sequel we  analyse the latter case in more
detail and identify the constant $\alpha$.

\begin{remark}
The relation $L=N^{-1}$ is appropriate if one considers the mean-field
model which  in the hierarchical mean-field limit is a potential-theoretic analogue of
two-dimensional space as $N \to \infty$, see \cite{DGsel} and for which
migration, mutation, selection are all {three together} important.
\end{remark}

The task is now to describe the system in the limit $N \to \infty$
in two time windows:
\begin{itemize}
\item times of order 1 after starting: $(X^N(t))_{t \geq 0}$
\item
times $\alpha^{-1} \log N$ after starting: $(X^N (\alpha^{-1} \log N+t)^+)_{t \in \R}$.
\end{itemize}
In the first window the evolution of the small advantageous droplet is of primary
interest, in other words the process $\gimel^{N,m}$,
while in the second case when the droplet covers a positive proportion of the whole space, the system is best described by the
empirical measure $\Xi_N$ (see \ref{angr5}) or the tagged sample see (\ref{angr5++}).

\subsection{The early time window as $N \to \infty$}
\label{ss.finitetw}

Here we want to describe (1) the evolution of the atomic measure-valued process
$\gimel^{N,m}$ in
the limit $N \to \infty$ over times in some finite interval and (2) the limit
$N \to \infty$ of the dynamic of the empirical measure.

{\bf (1)} {\em Droplet evolution }

During  the early times
there will be a finite random number of sites where the advantageous type has
mass exceeding some prescribed $\ve >0$. Therefore as $N \to \infty$ we
expect a limiting evolution of $(\gimel^N_t)_{t \geq 0}$.
First of all we can show

\beP{P.lt}{(Limiting droplet dynamic)}

As $N\to\infty$ \be{agrev2}
\CL [(\gimel^{N,m}_t)_{t \geq 0}]
\Nto \CL[({\gimel^m_t})_{t \geq 0}],
\ee in the sense of convergence of
continuous $\mathcal{M}_a([0,1])$-valued processes where the set of atomic measures
$\mathcal{M}_a([0,1])$ is equipped with the weak atomic topology.
The topology is introduced in \cite{EK4} and we don't touch this
further here.
\hfill$\qquad \square$
\end{proposition}

To identify the limit
evolution $\gimel^m$ we need a bit of classical excursion theory (compare \cite{PY}, section 3, \cite{Hu}).

\beL{L.Pa2} {(Single site: entrance and excursion laws)}

(a)
 Let $c>0, d>0,\; s>0$.  Then $0$ is an exit boundary for the
the Fisher-Wright diffusion
\be{onedim}
dx(t)=-c x(t)dt
+sx(t)(1-x(t))dt +\sqrt{d \cdot x(t)(1-x(t))}dw(t),
\ee
which then has a $\sigma$-finite entrance law from state $0$
at time 0, the $\sigma$-finite excursion law
\be{dd300}
{ \mathbb{Q}} = \mathbb{Q}^{c,d,s}\ee on
\be{dd30} W_0:=\{w\in C([0,\infty),\mathbb{R}^+),\;w(0)=0,\;
w(t)>0\mbox{ for  } 0<t<\zeta \mbox{ for some }
\zeta\in(0,\infty)\}.\ee

(b) Moreover, denoting by $P^\ve$ the law of the process started
with $w(0)=\ve$ and $\ve> 0, \Q$ is given by:
\be{dd31a}
 \mathbb{Q}(\cdot) =\lim_{\ve\to 0}\frac{P^\ve(\cdot)}{S(\ve)},
 \ee
where $S(\cdot)$ is the scale function of the diffusion
(\ref{onedim}), defined by the relation,
\be{dd31ab}
P_\ve(T_\eta<\infty)=\frac{S(\ve)}{S(\eta)},\qquad
0<\ve<\eta<\infty,
\ee
where $T_\eta$ is the first hitting time of $\eta$.

For the Fisher-Wright diffusion $S$ is given by (cf. \cite{RW},
V28)) the initial value problem:
\be{dd31b} S(0)=0,
\frac{dS}{dx}=\frac{e^{-2sx}}{(1-x)^{2c}}, \ee so that \be{dd31c}
\lim_{\ve\to 0}\frac{S(\ve)}{\ve}=1. \ee

(c) The measure $\Q$ is $\sigma$-finite, namely
for any $\eta >0, \zeta$ as in (\ref{dd30}),
\be{dd32} \mathbb{Q} (\{w:\zeta(w)>\eta\})<
\infty,\ee
\be{dd32b2}
\mathbb{Q}(\sup_t(w(t))>\eta)=\mathbb{Q}(T_\eta
<\infty)= \frac {1}{S(\eta)} \la \infty \mbox{ as } \eta \to 0,
\ee and \be{dd32b} \int_0^1 x
\mathbb{Q}(\sup_t(w(t))\in dx)=\infty.\qquad \square \ee
\end{lemma}

Now the limit $(N \to \infty)$ droplet dynamic $\gimel^m$ can be identified as follows.

\beP{CSB}{(A continuous atomic-measure-valued Markov process)}

Let
${ N(ds,da,du,dw)}$ be a Poisson random measure on (recall
(\ref{dd30}) for $W_0$) \be{agrev81} [0,\infty)\times [0,1]\times
[0,\infty)\times W_0, \ee with intensity measure
\be{dd32c}
 ds\, da\, du\,\mathbb{Q}(dw),
 \ee
where $\mathbb{Q}$ is the single site excursion law defined in (\ref{dd31a}) in Lemma
\ref{L.Pa2}.

Then the following two properties hold.

(a) The stochastic integral equation for $(\gimel^m_t)_{t \geq 0}$ given as
\be{ZL2m}
\gimel^m_t=\int_0^t\int_{[0,1]}\int_0^{q(s,a)}\int_{W_0}w(t-s)\delta_a
{ N(ds,da,du,dw)}, \quad t \geq 0, \ee where  $q(s,a)$ denotes the
non-negative predictable function
\be{dd32d}
q(s,a):= (m+ c
\gimel^m_{s-}([0,1])), \ee has a unique continuous
$\mathcal{M}_a([0,1])$-valued solution, which equals
$(\gimel^{m}_t)_{t \geq 0}$ from equation  (\ref{agrev2}).

(b) The process $({\gimel^m_t})_{t \geq 0}$ has the following properties:
\begin{itemize}
\item
branching property (the process is like a branching process with immigration)
\item the mass of each atom observed from the time of its creation
follows an excursion from zero generated from the excursion law
$\mathbb{Q}$ (see (\ref{dd31a})), \item new excursions are produced
at time $t$ at rate
\be{agrev3}
m+ c\gimel^m_{t \cdot}([0,1]), \ee
\item
each new excursion produces an atom located at a point  $a \in
[0,1]$ chosen according to the uniform distribution on $[0,1]$,
\item at each $t$ and $\ve >0$ there are at most finitely many
atoms of size $\geq \ve$.
$\qquad \square$
\end{itemize}
\end{proposition}

{\bf (2)} {\em McKean-Vlasov equation {for limiting empirical measure}.}

Turn now to the global description of the complete population by
the empirical measure.
In the limit $N\to\infty$ the evolution of the empirical measure in a finite time window is given
by the McKean-Vlasov limit of our system, but of course it is trivial, that is, totally concentrated on type $1$ if the given initial state has this property.
This is however different at late times.
Consider therefore the above system (\ref{angr1})-(\ref{angr3})
of $N$ interacting sites with type space $\mathbb{K}=\{1,2\}$ starting at
time $t=0$ from a product measure (that is, i.i.d. initial values
at the $N$ sites). The {\em basic McKean-Vlasov limit} (cf. \cite{DG99},
Theorem 9) says that if we start initially in an i.i.d.
distribution, then
\be{angr6} \{\Xi_N(t)\}_{0\leq t\leq T} \Nto
\{\mathcal{L}_t\}_{0\leq t\leq T}, \ee
where the
$\CP(\CP(\{1,2\})$-valued {\em deterministic} path $\{\mathcal{L}_t\}_{0\leq t
\leq T}$ is the law of a {\em nonlinear} Markov process solving a {\em forward}
equation, namely the unique weak
solution of the {\em McKean-Vlasov equation}:
\be{angr7b} \frac {d
\mathcal{L}_t}{dt}= ({L}_t^{\mathcal{L}_t})^\ast \mathcal{L}_t,
\ee where for $\pi \in \mathcal{P}(\mathcal{P}(\mathbb{K}))$,
$L^{\pi}$ is given by the generator of the process given by the
evolution of type 1 in  (\ref{angr7b0b}) below,
(and $\pi =m(t)$) and the $\ast$ indicates the adjoint
of an operator mapping from a dense subspace of
$C_b(E, \R)$ into $C_b(E,\R)$
w.r.t. the pairing of $\CP(E)$ and $C_b(E,\R)$
given by the integral of the function with respect to the measure.

As pointed out above in (\ref{angr5+}), in the special case
of the type set $\{1,2\}$, we can simplify by  considering the
frequency of type 2 only and by reformulating (\ref{angr7b})
living on $\CP (\CP(\{1,2\}))$ in
terms of $\mathcal{L}_t(2) \in \mathcal{P}[0,1]$.  This we carry
out now.

Namely we note that given the mean-curve \be{angr7b0a} m(t) =
\intl_{[0,1]} y\,\mathcal{L}_t(2)(dy), \ee the process
$(\CL_t(2))_{t \geq 0}$ is the {\em law} of the  solution of
(i.e. the unique weak solution) the SDE:
\be{angr7b0b} dy(t) = c(m(t)-y(t)) dt + sy(t)(1-y(t))dt +
\sqrt{dy(t)(1-y(t)} dw(t),
\ee
with $w$ being standard BM.
Then informally $(\CL_t)_{t \geq 0}$ corresponds to the
solution of the nonlinear diffusion equation. Namely for $t >0, $
$\CL_t(2)(\cdot)$ is absolutely continuous and for
\be{angr7b1}
 \mathcal{L}_t(2)(dx)=u(t,x)dx \;\in \mathcal{P}([0,1])
 \ee
 the evolution equation of the density $u(t,\cdot)$ is given by:
 \be{angr7b0}
  \frac{\partial}{\partial t} u(t,x) =
    - c\frac{\partial}{\partial x}\{[\intl_{[0,1]} y u(t,y) dy - x]{u(t,x)}\}
    -s\frac{\partial}{\partial x}(x(1-x)u(t,x))
      + \frac{d}{2}\frac{\partial^2}{\partial x^2}(x(1-x)u(t,x)).
      \ee

We have the following basic property for the McKean-Vlasov equation.

\beP{P.MCV}{(McKean-Vlasov equation and its solution)}

(a) Given the initial state $\mu_0 \in \CP([0,1])$
there exists a unique solution
\be{angr8b} \mathcal{L}_t(2)(dx) = \mu_t(dx),\;\; t\geq t_0 \ee to
(\ref{angr7b}) with initial condition $\mathcal{L}_{t_0}(2)=\mu_0$. \\
(b) If $s>0$ and $ \int_{[0,1]} x\mu_{t_0}(dx) >0$, then this
solution satisfies:
\be{angr8c}
   \lim_{t\to\infty}
   \mathcal{L}_t(2)(dx) =\delta_1(dx).\qquad \square
   \ee
\end{proposition}

\subsection{The late time window as $N \to \infty$}
\label{ss.latetw}

In the late time window we see (global) {\em emergence} and then {\em fixation}
of the advantageous type.

{{\bf (1)} Emergence times.}
Here we begin by studying the emergence by identifying
the time of emergence and of fixation of the advantageous type as
follows.
\beP{P6.1b}{(Macroscopic emergence and fixation times)}

\noi(a) {(Emergence-time)}

There exists a constant $\alpha$ with:
\be{ang3}
0<\alpha<s, \ee
such that if
$\; T_N = \frac{1}{\alpha}\log N$, then for $t \in \R$ and asymptotically as
$N \to \infty$ type 2 is present at times $T_N+t$, i.e. there exists a $\ve >0$
such that for every $i$,
\be{Y12b-}
  \liminf_{N \to \infty}  P[ x^N_2 (i, T_N+t)>\ve]>0,  \ee
  and type 2 is not present earlier, namely for   $1>\ve >0$:
  \be{Y12b}
   \lim_{t\to -\infty} \limsupl_{N \to \infty}
    [P( x^N_1 (i, T_N+t)<1-\ve]=0.  \ee

\noi(b) (Fixation time)

After emergence the fixation occurs in times $O(1)$ as $N \to \infty$,
i.e. for any $\ve>0, \quad i \in \N$
\be{Y12c}
\lim_{t\to\infty}\limsup_{N\to\infty}P[x^N_1(i, T_N+t)>\ve]=0. \qquad \square
\ee
\end{proposition}

\beC{C.emfix}{(Emergence and fixation times of spatial density)}

The relations (\ref{Y12b-}), (\ref{Y12b}) and (\ref{Y12c})
hold for $\bar x^N_2$ respectively $\bar x^N_1$ as well.
$\qquad \square$
\end{corollary}

We can identify the parameter $\alpha$ as follows from the droplet growth
behaviour.

\beP{CSB-longtime}{(Long-time growth behavior of $\gimel^m_t$)}

Assume that  $m>0$.  Then the following
growth behavior of $\gimel^m$ holds.\\
(a) There exists $\alpha^\ast$ such that the following limit exists
\be{CSB1}
\lim_{t\to\infty} e^{-\alpha^\ast t}E[\gimel^m_t([0,1])] \in (0,\infty),
\ee
with { (here $\alpha$ is from (\ref{ang3}))}
\be{CSB1b}
\alpha^\ast = \alpha, \text{ where }\alpha \text{ is given below in (\ref{ang4c4})}
\ee

\be{CSB3} e^{-\alpha t}
\gimel^m_t([0,1]) \tto {\mathcal{W}}^\ast \mbox{ in probability},
\; 0<\CW^\ast < \infty \mbox{ a.s.}.
\ee

(b) The growth factor in the exponential is truly random:
\be{CSB2}
0< Var[\lim_{t\to\infty} e^{-\alpha t}\gimel^m_t([0,1])]<\infty. \qquad \square
\ee

\end{proposition}

\begin{remark}
The random variable $\CW^\ast$
reflects the growth of $\gimel^m_t([0,1])$ in the beginning, as
is the case in a supercritical branching process and hence
$\CE^\ast = \alpha^{-1} \log \CW^\ast$ can be viewed as the random
time shift of that exponential $e^{\alpha t}$ which matches
the total mass of $\gimel^m_t$ for large $t$.
\end{remark}

\begin{remark}
Even though one might think that for small mass of the advantageous type this
expands at the rate $e^{st}$, this is not the case due to the stochastic effects
leading to a subtle interplay between the parameters $s,d$ and $c$ resulting
in $\alpha^\ast = \alpha <s$.
\end{remark}

{{\bf (2)} {\em Fixation dynamic.}}
We now understand the preemergence situation
and the time of emergence and fixation.
In order to describe the whole dynamics of
{\em macroscopic fixation}, we  consider the limiting
distributions of the empirical measure-valued processes in a
second time window $(\alpha^{-1} \log N+t)^+, \quad t \in \R$. Define
(we suppress the truncation in the notation below)
\be{angr9z}
\Xi^{\log,\alpha}_N(t):=\frac{1}{N}
\sum_{i=1}^N \delta_{(x_1^N(i,\frac{\log
N}{\alpha}+t),x_2^N(i,\frac{\log N}{\alpha}+t))},\quad t \in
\R,\quad(\Xi^{\log,\alpha}_N(t)\in \CP(\mathcal{P}(\mathbb{K}))) \ee
and then the two empirical marginals are given as
\be{angr9}
\Xi^{\log,\alpha}_N(t,\ell):=\frac{1}{N} \sum_{i=1}^N
\delta_{x_\ell^N(i,\frac{\log N}{\alpha}+t)}, \; \ell=1,2 \mbox{ and } t
\in \R, \quad(\Xi^{\log, \alpha}_N(t,\ell)\in \mathcal{P}({[0,1]})).
\ee Note that for each $t$ and given $\ell$ the latter is a random
measure on $[0,1]$. Furthermore we have the representation of the
empirical mean of type 2 as follows:
\be{angr10}
\bar x^N_2(\frac{\log
N}{\alpha}+t) =\int_{[0,1]} x\;\Xi^{\log,\alpha}_N(t,2)(dx).\ee

Since we consider the limits of systems observed in the interval
$const\cdot \log N +[-\frac{T}{2},\frac{T}{2}]$ with $T$ any
positive number, that is setting $t_0(N) = const \cdot \log N
-\frac{T}{2}$,  we need to identify entrance laws for the process
from $-\infty$ (by considering $T \to \infty$) out of the state
concentrated on type 1 with certain additional properties.

The next main result is on the fixation process,
saying that $\Xi^{\log,\alpha}_N$
converges as $N\to\infty$ and that the limit can be explicitly
identified as a random McKean-Vlasov entrance law from $-\infty$,
a concept we explain next.

\beD{D.entlaw} {(Entrance law from $t = -\infty$)}

We say in the two-type case that a probability measure-valued function
$\CL : \R \to \CP([0,1])$, is an entrance law at $-\infty$
starting from type 1 if
$(\CL_t)_{t \in (-\infty,\infty)}$ is
such that $\CL_t$ solves the McKean-Vlasov equation (\ref{angr7b}) and
$\CL_t \to \delta_1$ as $t \to -\infty$.
 $\qquad \square$
 \end{definition}

We will indeed establish that the emergence of rare mutants gives
rise to \emph{``random''} solutions of the McKean-Vlasov dynamics.
In particular we will show that the limiting empirical measures at
times of the form $C \log N +t$ are {\em random} probability measures on
$[0,1]$ and therefore given by sequences of $[0,1]$-
valued truly exchangeable random variables which are \emph{not} i.i.d., that is,
the exchangeable $\sigma$-algebra is not trivial.
This means that the limiting empirical mean turns out to be a {\em random variable}
and this is the  driving term due to migration for the local evolution
of a site in the McKean-Vlasov limit.  Therefore both the random driving term and the non-linearity of the evolution
equation come seriously into play. However once we condition on the
exchangeable $\sigma$-algebra, we then get for the further evolution
again a deterministic limiting equation for the empirical measures,
namely the McKean-Vlasov equation. The reason for this is the fact
that conditioned on the exchangeable $\sigma$-algebra we obtain an
asymptotically (as $N \to \infty$) i.i.d. configuration to which
the classical convergence theorem applies. Using the Feller property
of the system, we get our claim.
This leads to
the task of identifying an entrance law in terms of a {\em random
initial condition at time $-\infty$}.

The above discussion shows that we need to
introduce  the notion of a truly random McKean-Vlasov entrance
law from $-\infty$.

\beD{D.ranMKV}{(Random entrance laws of McKean-Vlasov from $t =-\infty$)}

We say that the  probability measure-valued process
$\{\CL(t)\}_{t\in\mathbb{R}}$ is a {\em random solution of the
McKean-Vlasov equation} (\ref{angr7b}) if
\begin{itemize}
\item $\{\CL_t:t\in \R\}$ is a.s. a solution to (\ref{angr7b}),
that is, for every $t_0$ the distribution of $\{\CL_t:t\geq t_0\}$
conditioned on $\mathcal{F}_{t_0}=\sigma\{\CL_s:s\leq t_0\}$ is
given by $\delta_{{\{\mu_t}\}_{t\geq t_0}}$ where $\mu_t$ is a
solution of the McKean-Vlasov equation with $\mu_{t_0}=\CL_{t_0}$,
\item the time $t$ marginal distributions of $\{\CL_t : t \in
\R\}$ are truly random.$\qquad \square$
\end{itemize}
\end{definition}
We can say the following about the possible random entrance laws.

\beP{P.random}{(Random entrance laws)}

(a)
 There exists a solution
 $(\mathcal{L}_t^{\ast \ast}(2))_{t \in \R}$ to  equation
 (\ref{angr7b})
 satisfying the conditions:
 \bea{Ex1b}
 \lim_{t\to-\infty} \mathcal{L}^{\ast\ast}_t(2) &= &\delta_0,\\
 \lim_{t\to\infty} \mathcal{L}^{\ast\ast}_t(2) &= &\delta_1 \nonumber \\
 \int_{[0,1]}x\mathcal{L}^{\ast\ast}_0 (2,dx) &=
 &\frac{1}{2}.\nonumber \eea This solution is called an entrance
 law from $-\infty$ with mean
 $\frac{1}{2}$ at $ t=0$.\\
 (b) We can obtain a solution in (a) such that:
 \be{Ex2}
 \exists\; \alpha
 \in (0,s)\mbox{ and  }A_0\in(0,\infty)\mbox{ such that }
 \lim_{t\to-\infty}e^{\alpha |t|}
  \int_{[0,1]} x \mathcal{L}^{\ast\ast}_t(2,dx) =A_0.
  \ee
  (c) The solution of (\ref{angr7b}) also satisfying
  (\ref{Ex2}) for prescribed $A_0$ is unique and if $A_0 \in (0, \infty)$
  then $\alpha$ is necessarily uniquely determined.
						
For any deterministic solution
\be{Ex2a}
\{\mathcal{L}_t,\; t \in \R\}
\ee
 to (\ref{angr7b}) with
 \be{Ex2b} 0 \leq \limsup_{t\to -\infty}e^{\alpha |t|}\int_{[0,1]}
 x \mathcal{L}_t(2,dx) <\infty, \ee the limit
 $A=\lim_{t\to-\infty}e^{\alpha |t|} \int_{[0,1]} x
 \mathcal{L}_t(2,dx)$ exists.
			
If $A>0$, then $\{\mathcal{L}_t, t \in \R\}$ is
given by a time shift of the then unique
$\{\mathcal{L}^{\ast \ast}_t, t \in \R\}$ singled out in (\ref{Ex2}), i.e.
\be{Ex2c}
\CL_t = \CL^{\ast\ast}_{t+\tau}, \quad \tau =
\alpha^{-1} \log \frac{A}{A_0}.
\ee
For future reference we define
$(\CL^\ast_t)_{t \in \R}$ to be the unique solution satisfying
\be{L.ast}\lim_{t\to-\infty}e^{\alpha |t|}
 \int_{[0,1]} x \mathcal{L}^\ast_t(2,dx)=1
 \mbox{ for some } \alpha \in (0,s).
 \ee

 (d) Any random solution $({\mathcal{L}}_t)_{t \in \R}$ to (\ref{angr7b})
 such  that
 \be{Ex2d}
 \limsup_{t\to -\infty}e^{\alpha |t|}E[\int_{[0,1]} x \mathcal{L}_t(2,dx)] <\infty,\quad
 \liminf_{t\to -\infty}e^{\alpha |t|}[\int_{[0,1]} x
 \mathcal{L}_t(2,dx)] >0 \mbox{ a.s.},
 \ee is a random time shift of
 $(\mathcal{L}^{\ast\ast}_t)_{t \in \R}$ and of $(\CL^\ast_t)_{t \in \R}$.
  \qad
\end{proposition}

\begin{example}  Let $\mathcal{L}^\ast_t$ be a solution
satisfying (\ref{angr7b0}), (\ref{Ex2}) and for a given value of
$A$  let $\tau$ be a true real-valued random variable.  Then
$\{\mathcal{L}^\ast_{t-\tau}\}_{t\in\mathbb{R}}$ is a truly random
solution. This can also be viewed as saying that we have a
solution with an exponential growth factor $A$ which is truly
random.
\end{example}

The emergence of the advantageous type and the subsequent evolution
to fixation in this type is characterized as follows.

\beP{P.measure} {(Asymptotic macroscopic emergence-fixation process)}

(a) For each $-\infty<t<\infty$ the empirical measures converge
weakly to a random measure:
\be{angr11}
\mathbf{\CL} [\{\Xi^{\log,\alpha}_N(t,\ell)\}] \Nto \mathbf{\CL}[\{ \CL_t(\ell)\}] =P^\ell_t \in
\CP(\CP({[0,1]})), \mbox{ for } \ell =1,2 . \ee
In addition we have path convergence:
\be{angr13}
 w-\liml_{N \to \infty}\mathbb{\CL}[(\Xi^{\log,\alpha}_N(t))_{t \in \R}]=P
  \in \CP [C((-\infty,\infty),\CP(\CP(\{0,1\})))].
  \ee
  A realization of $P$ is denoted $(\CL_t)_{t \in \R}$
  respectively its marginal processes $(\CL_t(1))_{t \in \R}, (\CL_t(2))_{t \in \R}$.
  \sm

(b) The process $(\CL_t)_{t \in \R}$ describes the emergence and
fixation dynamics, that is, for $t\in\mathbb{R}$, and $\ve
>0$, \be{LD2} \lim_{t\to-\infty}
\mbox{Prob}[\mathcal{L}_t(2)((\ve,1])>\ve]=0,\ee
\be{LD3}
\lim_{t\to\infty}
\mbox{Prob}[\mathcal{L}_t(2)([1-\ve,1])<1-\ve]=0,\ee
with
\be{LD3b}
\CL_t (2) ((0,1)) > 0\quad , \quad \forall t \in \R, \mbox{ a.s.}.
\ee
(c) The limiting dynamic in (\ref{angr13}) is identified as follows:

The probability measure $P$ in (\ref{angr13}) is such
that the canonical process is a  random solution  (recall
Definition \ref{D.ranMKV}) and entrance law from time $-\infty$ to the McKean-Vlasov equation
(\ref{angr7b}). \sm

(d) The limiting dynamic in (\ref{angr13}) satisfies with $\alpha$ as in (\ref{ang3}):
\be{angr14a2}
\CL [e^{\alpha |t|}\int_{[0,1]} x \CL_t(2)
(dx)]\Rightarrow \CL[^\ast \CW] \mbox{ as } t \to -\infty,
\ee
 and we explicitly identify the
 random element generating $P$ in (\ref{angr13}), namely $P$ arises
 from random shift of a deterministic path:
 \be{angr14c}
 P = \CL [\tau_{^\ast\CE} \CL^\ast]
  \quad , \quad ^\ast\CE =(\log ^\ast\CW)/\alpha, \quad
   \tau_r \mbox{ is the time-shift of path by } r,
   \ee
   where $\CL^\ast$ is the unique and deterministic entrance law of the McKean-Vlasov
   equation (\ref{angr7b})  with projection $\CL^\ast_t(2)$ on the type 2
   coordinate satisfying:
   \be{angr14d} e^{\alpha |t|} \int_{[0,1]} x
   \CL^\ast_t(2)(dx) \longrightarrow 1,  \mbox{ as } t \to -\infty.
    \ee

 The random variable $^\ast\CW$ satisfies \be{angr14a}
 0< {^\ast{\CW}} < \infty \mbox{ a.s.},\quad E[^\ast\CW]<\infty,\quad
 0<\mbox{Var}(^\ast\CW)<\infty. \ee

(e) We have for $s_N \to \infty$ with $s_N =o(\log N)$ the approximation
    property for the growth behaviour of the limit dynamic by the finite $N$ model,
    namely:
    \be{angr14a3}
    \CL [e^{\alpha s_N} \bar x^N_2 (\frac{\log N}{\alpha} - s_N)]
    \Nto \CL[^\ast\CW]. \qquad \square \ee
    \end{proposition}

How does this emergence in (\ref{angr14a3}) relate to the droplet growth? We have

\beP{P.DF}{(Microscopic emergence and evolution: droplet formation)}

The total type-2 mass $\gimel^{N,m}_{t_N} ([0,1])$ grows at exponential rate
$\alpha$,
\be{rgf1}
\CL \left[\gimel^{N,m}_{t_N} ([0,1])
e^{-\alpha t_N}\right] \Nto \CL [ \CW^\ast], \quad \mbox{ for }
t_N \uparrow \infty \mbox{ with } t_N -(\alpha^{-1} \log N) \la -\infty,
\ee
\be{ad7}
\CL[^\ast\CW] = \CL[\CW^\ast]. \qquad \square
\ee
\end{proposition}

\section{A logistic branching random walk and its growth }
\label{s.logisticbrw}

We define here a logistic branching population, take the mean-field limit
$N \to \infty$, study the expansion of the droplet of occupied sites, determine
the late time window and finally
the behaviour in a late time window (i.e. shift of  observation time interval and
size of space tend to infinity) in four subsections.

\subsection{The logistic branching particle model}
\label{ss.partmod}

We now consider a particle system on the geographic space
$\{1,\cdots,N\}$ and its occupation number configuration with state space
\be{ad4}
(\N_0)^{\{1, \cdots,N\}}.
\ee
Here particles
\begin{itemize}
\item
{\em migrate} according to a continuous time (rate $c$) random walk
with uniform step distribution,
\item
a particular particle at a site $i$ {\em dies} at rate $d(k-1)$ if we have $k$
particles at site $i$,
\item
a particle has {\em one offspring} at rate $s$ which is placed at the same site.
\end{itemize}

This defines uniquely a strong Markov pure jump process on our state space, which
we denote by
\be{ad5}
(\eta^N_t)_{t \geq 0}.
\ee
\begin{remark}
Note that the mean production rate (mean growth rate) of this model in state $k$ is
\be{ad6}
sk - d(k(k-1)),
\ee
which is a concave function $f$ with $f(0)=0, f(1)=s 1, f(k) >0$
for $k \leq k_0$ and $f(k)<0$ for $k >k_0$.
This production rate can be interpreted as reflecting limited local resources for a
population.
\end{remark}

\begin{remark}
We can view this process as a supercritical branching random walk (supercriticality parameter $s$)
with  an additional linear  death rate $d(k-1)$ per individual thus inducing an interaction between families.
The process is also called a coalescing branching random walk.
\end{remark}

Due to the quadratic death rate the population can only expand indefinitely by having  individuals
move to sites where so far no particles are sitting. When the space is filled
the expansion of the population is replaced by an equilibrium situation. How
to make this precise?

We want to study this particle system for $N \to \infty$ in finite time
windows, one early starting at time 0 and the other time windows beginning at a late
time {when space fills up with particles}.
Here we start with one particle and determine the late time window
by asking when does the population develop a positive spatial intensity
even as $N \to \infty$.
Then in particular we focus on the influence of the deathrate $d$ and
the migration rate $c$ on the speed of spatial spread.
The goal then is to establish that  the population grows
exponentially fast with a rate $\alpha$ which is positive
but strictly less than the birth rate $s$ as long as we have no collisions.

\begin{remark}
Note that if $d=0$, then in fact we have a supercritical branching process
at exponential rate $s$ and it would take time
$\frac{1}{s}\log N$ to develop positive spatial intensity.
\end{remark}

If we have only one site we have a classical birth and death process,
in the spatial model we have a collection of such processes.
Since the death rate is zero if we have only one particle at a site, and furthermore
because the  death rate is quadratic, we have a positive recurrent
Markov process on the state space $\N^N_0$ with a unique equilibrium law denoted:
\be{ad8}
\pi^N = \pi^N_{c,s,d}. \mbox{ Write }
\pi_{c,s,d} \mbox{ for } \pi^1_{c,s,d}.
\ee
Out of the finite initial  state a new site can be occupied and the population can
grow till the jumps can only hit already occupied sites and this way a
population intensity on the whole space can develop and a local equilibrium
forms, resulting in a global equilibrium density.

\subsection{The early time window as $N \to \infty$}
\label{ss.ftwinfty}

As $N \to \infty$ the geographic space expands to $\N$ and the migrating particles
eventually do not  hit occupied sites in a fixed time interval (if we start with finitely many
initial particles, the collision-free regime, later when collisions occur sites
interact again). More precisely,
we want to establish that, as $N \to \infty$, we get as limit dynamic a collection
of birth and death processes with emigration at rate $c$ and immigration
at rate $c$ from a reservoir with the current intensity in the total
population. This is carried out  as follows.

Consider a branching random walk with birth and death as
before but now with geographic space $\N$, so that the state space
becomes
\be{ad9}
(\N_0)^N
\ee
and migration changes to
\begin{itemize}
\item
emigration out of any component at rate $c$ to the unoccupied site
of lowest index as long as there is more than one particle
\item
immigration at rate $c \cdot \iota$ into every colony, with $\iota \in [0,\infty)$.
\end{itemize}
Here if we start the system in an exchangeable initial state we choose
\be{ad10}
\iota = \lim_{N \to \infty} (N^{-1} \suml^N_{i=1} \eta (i)).
\ee
Later in larger time scales we will obtain a $\iota$ which is time-dependent
and arises from the law at a tagged site. For $\iota=0$ we obtain
what we call the {\em collision-free process}.

The strong Markov process defined by the above McKean-Vlasov dynamic is denoted
\be{ad11}
(\eta^\iota_t)_{t \geq 0} \quad (\mbox{ resp. } (\eta_t)_{t \geq 0} \mbox{ if } \iota=0).
\ee
This process has for every value $\iota \in [0,\infty)$ a unique
equilibrium
\be{ad12}
\Pi^\iota = \bigotimes \pi^{(\iota)}_{c,s,d},
\ee
where $\pi^\iota_{c,s,d}$ is the single site equilibrium.
Let $r$ be any map which permutes the location such that
$1,2,\cdots, |\{i|\eta^N_t(i)>0\}|$
are all occupied {and acts on paths by achieving the constraint at the
final time}.  We prove
\beP{P.convcfp}{(Convergence to $\iota$-process, collision-free process)}

If we start in an i.i.d. distribution with
$E[\eta^N_t(i)] = \iota$, then
\be{ad26}
\CL [(\eta^N_t)_{t \geq 0}] \Nto \CL [(\eta^{(\iota)}_t)_{t \geq 0}].
\ee
If we start with one initial particle at site 1 and if
$t_N=o(\log N)$, then
\be{ad27}
\CL [(r \circ \eta^N_t)_{t \geq 0}] \Nto \CL [(\eta_t)_{t \geq 0}].
$\qad$
\ee
\end{proposition}

\subsection{The droplet expansion and Crump-Mode-Jagers processes}
\label{ss.dropletexp}

As time gets large we expect the populations of particles to occupy
more and more sites (droplet), such that the overall population expands as time
grows, i.e. we have a growing {\em droplet of occupied sites}. This droplet growth
is expected to be exponential and as time grows to become (at least
the further evolution) more and more deterministic due to a law
of large numbers effect. Here we have to distinguish a time window,
which is late but where particles essentially always move to new
unoccupied sites if they migrate and a later phase where spatial
intensity builds up and collisions play a role. This subsection
handles the first regime.

If we establish an exponential growth behaviour we would be able to determine the
position of the late time window where the population density
becomes positive as time and $N$ tend to infinity.
How can we make all this rigorous?

To analyse this time window an important concept is that of a CMJ-process
(Crump-Mode-Jagers process) which allows us to describe the number of
occupied sites in the process $(\eta_t)_{t \geq 0}$ starting from
finitely many particles occupying all sites $1, \cdots, k$ with $k \geq 1$. Let
\be{ad13}
K_t = \# \{i \in \N | \eta_t(i) >0\},
\ee
\be{ad14}
\zeta_t = \{\zeta_t (i), \quad i \in 1,\cdots, K_t\} \quad ,\quad \zeta_t(i)=\eta_t(i),
\quad i=1, \cdots, K_t.
\ee
We will find that $K, \zeta$ are processes which have the structure of
a CMJ-process that is the process of {\em occupied sites} is a type of
generalised branching process. We begin by recalling this concept.
\sm

{\bf 1.}
{\em Crump-Mode-Jagers process}

The CMJ-process models individuals in a branching
population whose dynamics is as follows.
Individuals can die or give birth to new individuals
based on the following ingredients:
\begin{itemize}
\item individuals have a lifetime (possibly infinite), \item for
each individual an independent realization of a point process
$\xi(t)$  starting at the birth time specifying the times at which
the individual gives birth to new individuals, \item different
individuals act independently, \item the process of birth  times is
not  concentrated on a lattice.
\end{itemize}

The process might be growing exponentially and then we want to
determine its exponential growth rate.
The corresponding  {\it Malthusian parameter}, $\alpha
>0 $ is obtained as the unique solution of
\be{ang1}
\int_0^\infty
e^{-\alpha t}  \mu (dt) =1 \mbox{ where } \mu([0,t]) =
E[\xi([0,t])],
\ee
with
\be{ang2}
\xi(t) \mbox{ counting the number of births of
a single individual up to time } t,
\ee
(see for example \cite{J92}, \cite{N} equation (1.4)).
\sm

If we know the Malthusian parameter, we need to know that it is
actually equal to the almost sure growth rate of the population.
It is known for a CMJ-process $(K_t)_{t \geq 0}$ that (Proposition 1.1
and Theorem 5.4 in \cite{N}) the following basic growth theorem
holds. If
\be{ang1b}
E[X \log (X \vee 1)] < \infty, \quad \mbox{ where }
X =\intl^\infty_0 e^{-\alpha t} d \xi(t),
\ee
then
\be{ang2b}
 \lim_{t\to\infty} \frac{K_t}{e^{\alpha t}} = W,  \mbox{ a.s. and in } L_1, \ee
where $W$ is a random variable which has two important properties, namely
\be{ang4}
W>0,\;a.s.$ and $E[W]<\infty. \ee

The empirical distribution of individuals of a certain age in the population converges
as $t \to \infty$  to a {\em stable age  distribution}

\be{ang4c}
\CU(\infty,du) \mbox{ on } [0,\infty), \ee according to
Corollary 6.4 in \cite{N}, if condition 6.1 therein holds.
The condition 6.1 in \cite{N}  or (3.1) in \cite{JN} requires
that (with $\mu$ as in (\ref{ang1})):
\be{ang4c1} \int_0^\infty e^{-\beta t}\mu(dt)<\infty \quad
\mbox{ for some } \beta < \alpha.
\ee

{\bf 2.}
{\em Application of the CMJ-theory to $(K_t, \zeta_t)$}

The $(\zeta_t)_{t \geq 0}$ we have introduced is a CMJ
where individuals are the occupied sites and satisfying all the
conditions posed above. In addition we have more structure, namely the
birth process is determined from an {\em internal state} of the
individual (here a site) which follows a Markovian evolution
(here our $(\zeta_t(i))_{t \geq 0})$. This allows to obtain some
stronger statements as follows.

Define
\be{ad15}
(\zeta_0(t,i))_{t \geq 0}, \mbox{ as the process }
(\zeta_t(i))_{t \geq 0}) \mbox{ if $i$ gets occupied at time } 0.
\ee
In our case we then have
 \be{ang1ab} \mu([0,t]) = c \int_0^t
E[\zeta_0(s) 1_{(\zeta_0(s) \geq 2)}]ds. \ee

Define the random measure
\be{ad16}
\CU(t, du,j) = \# \{\mbox{ sites of size $j$ with birth time } t-du\} K^{-1}_t.
\ee
It follows that the random measure $ \CU(t,\cdot,\cdot)$ converges to a
{\em deterministic} object, the {\em stable  age
and size} distribution, i.e.
\be{ang4c2}
\CU (t,\cdot,\cdot) \La \CU(\infty,\cdot,\cdot), \mbox{ as } t\to \infty \mbox{ in law},
\quad \CU(\infty, \cdot, \N) = \CU(\infty, \cdot) \mbox{ from (\ref{ang4c})}.
\ee
We can obtain from this the following representation of the
Malthusian parameter $\alpha$.

Given site $i$ let
$\tau_{i}\geq 0$ denote the time at which a migrant (or initial
particle)  first occupies it. Noting that we can verify Condition
5.1 in \cite{N} we have that
\be{stableage}
\lim_{t\to
\infty}\frac{1}{ K_t}\sum_{i=1}^{ K_t} \zeta_{\tau_i}(t-\tau_i)
=\intl^\infty_0 E [\zeta_0(u)] \CU(\infty, du)=  B \text{  (a
constant)}, a.s.,
\end{equation} by Corollary 5.5 of
\cite{N}. The constant $B$ in (\ref{stableage}) is in our case given by the
average number of particles per occupied site and the growth rate $\alpha$
arises from this quantity neglecting single occupation. Namely define

\be{ang4c4}
\alpha =c\sum_{j=2}^\infty j
\mathcal{U}(\infty,[0,\infty),j)<\infty ,\qquad
\gamma=c \; \mathcal{U}(\infty,[0,\infty),1).\ee
Then  the mean occupation
number of an occupied site is given by

\be{ang4c3}
B = \frac{\alpha+\gamma}{c} >0.
\ee
Furthermore  the average birth rate of new sites (by arrival of a
migrant at an unoccupied site) at time $t$ (in the process
in the McKean-Vlasov dual) is equal to
\be{ang4c5}
\alpha = c \; B - \gamma.
\ee

Here a remarkable point is the relation between $\alpha$ and $s$.
Recall $s$ is the parameter of supercriticality in the branching part of the internal state dynamics.
The action of the quadratic death part and the role of migration lead
to a number $\alpha$ with
\be{ad17}
0< \alpha <s,
\ee
since the fraction of the supercritical branching process with
supercriticality parameter $s$ which die by the quadratic death rate is positive.

\subsection{Time point of emergence as $N \to \infty$}
\label{ss.timepoint}

If we take the collision-free model
$(\eta_t)_{t \geq 0}$ and we observe at time
$T_N(t) = \alpha^{-1} \log N +t$ the number of individuals
or the number of occupied sites, both normalized by $N$,
then these quantities satisfy (recall (\ref{ang4c3}), (\ref{ang2b}))
\be{ad18}
\left(\frac{1}{N} \suml^{K_{T_N(t)}}_{i=1} \eta_{T_N(t)} (i)\right)
\Ntoo B  e^{\alpha( t +\alpha^{-1} \log W)} ,
\quad (K_{T_N(t)}/N)_{t \in \R} \Ntoo
(e^{\alpha(t+\alpha^{-1} \log W)})_{t \in \R},
\ee
with $0<W<\infty$ a.s. and $Var (W) >0$. Hence a positive intensity develops:
\be{ad19}
\left( \frac{1}{K_{T_N(t)}} \suml^{K_{T_N(t)}}_{i=1} \eta_{T_N(t)}(i)\right)
\Ntoo B, \quad \forall \; t \in \R.
\ee
In particular in this time window $\eta^N$ and $\eta$ differ
significantly, but one can show that still:
\be{ad29}
\CL [(\{r \circ \eta^N_{T_N+t} (1), \cdots, ( r \circ \eta^N_{T_N+t} (k)\})_{t\geq 0}]
\La \CL[(\eta_t(1), \cdots, \eta_t(k))_{t \geq 0}], \mbox{ as } N \to \infty,
\ee
provided $T_N - \alpha^{-1} \log  N \to -\infty$
as $N \to \infty$. In fact even though the CMJ-approximation breaks down
at times $\alpha \log N+t$ we can still prove emergence occurs at time
$\alpha^{-1} \log N$, since $\eta^N$ develops
a positive intensity meaning the total population is comparable to $N$,
which as $t \to -\infty$ is asymptotically equivalent to the r.h.s. of
(\ref{ad18}).
Indeed time $\alpha^{-1} \log N$ separates the collision-free
droplet growth from the emergence.

\beP{P.empos} {(Emergence of positive intensity)}

We have for $N \to \infty$ and $\nu = \nu((t_N)_{N \in \N})$:
\be{ad25}
\CL[N^{-1} \suml^N_{i=1} \eta^N_{t_N} (i)]
\tto \left\{\begin{array}{l}
\delta_0 \quad , t_N - \alpha^{-1} \log N \to -\infty\\
\nu, \nu((0,\infty))=1, \quad \liml_{N\to \infty} (t_N- \alpha^{-1} \log N)=t > -\infty.
\end{array}
\right. \qquad \square
\ee
\end{proposition}

\subsection{The late time window as $N \to \infty$}
\label{ss.ltwinfty}

To treat the time window $ \alpha^{-1} \log  N+t$
and to obtain the limiting emergence-equilibration dynamics
based on {\em random entrance laws} of the {\em McKean-Vlasov equation}, we
need some ingredients.

We consider the number of sites occupied at time $t$ denoted $
K^N_t$ and the corresponding measure-valued process on
$[0,\infty)\times \N_0$ giving the unnormalized
number of sites of a certain age $u$ in $[a,b)$ and occupation size $j$:
\be{ba10d2}
\Psi^N(t,[a,b),j)) = \int_{(t-b)}^{(t-a)}1_{( K^N_{u}>
K^N_{u-})}1_{(\zeta^N_u(t)=j)}d K^N_u,
\ee
where $\zeta^N_u(t)$
denotes the occupation number at time $t$ of the site born at time
$u$, that is, a site first occupied the last time at time $u$, which is
therefore at time $t$ exactly of age $t-u$.

The {\em normalized empirical age and size distribution} among the occupied sites
is defined as: \be{ba10f1} U^N(t,[a,b),j) =\frac{1}{K^N_t} \Psi^N
(t,[a,b),j),\quad t\geq0,\; j\in \{1,2,3,\dots\}. \ee

\noi
Denote now for convenience  the  number of sites:

\be{Grev39c2}
 u^N(t):= K^N_t.
\ee

We have obtained with (\ref{ba10f1}), (\ref{Grev39c2}) a pair which is
$\N \times \CP([0,\infty) \times \N)$-valued
and which describes  our particle system completely provided
if we consider individuals and sites as exchangeable labels.
This pair is denoted
\be{ba10f4}
(u^N(t), U^N(t,\cdot,\cdot))_{t \geq 0}.
\ee

In the interval $\alpha^{-1} [\log N- \log \log N,
\log N+T]$ the process $(u^N(t))_{t \geq 0}$  increases by one
respectively decreases by one at rates

\be{Grev39d}
 \alpha_N (t)(1-\frac{u^N(t)}{N}) u^N(t), \mbox{ respectively }
\gamma_N (t) \frac{( u^N(t))^2}{N}, \ee
where $\alpha_N
(t),\gamma_N (t)$ are defined: \be{Grev39e} \alpha_N (t) = c
\intl^t_0 \suml^\infty_{j=2} jU^N (t,ds,j) ds, \quad \gamma_N (t)
= c \intl^t_0 U^N (t,ds,1). \ee
The above rates of change of $U^N (t,\cdot,\cdot)$ follow directly
from the dynamics of the  particle system $\eta^N$.

Our goal is now to study the behaviour of
$(u^N, U^N)$ at times $\alpha^{-1} \log N+t$ and to show that
this follows a limiting fixation dynamics. We start the system with
$k$ particles at $\ell$ distinct sites and write $k, \ell$ as superscript. We need
the time-shifted quantities:
\be{ts1}
 \wt u^{N,k,\ell}(t)
  = u^{N,k,\ell}((\frac{\log N}{\alpha}+t)\vee 0),\; t\in (-\infty,T],\quad \wt u^{N,k,\ell}
  (-\frac{\log N}{\alpha})=\ell\ee
  \be{ts2}
  \wt U^{N,k,\ell}(t)
   = U^{N,k,\ell}((\frac{\log N}{\alpha}+t)\vee 0),\; t\in (-\infty,T],\quad \wt U^{N,k,\ell}
   (-\frac{\log N}{\alpha})=
    \delta_{(k,0)}.\ee

\beP{P.Grocoll}{(Convergence to a colonization-equilibration dynamic in the $N\to\infty$ limit)}

Assume that for some $t_0 \in \R$ as $N\to\infty$, $(\frac{1}{N} \wt u^{N,k,\ell}
(t_0), \wt U^{N,k,\ell} (t_0))$ converges in law to the pair
$(u(t_0),U(t_0))$ automatically contained in $[0,\infty) \times L_1(\N, \nu)$.

Then as $N \to \infty$
\be{agrev60}
\CL \left[\big(\frac{1}{N} \wt u^{N,k,\ell}(t),\,
 \wt U^{N,k,\ell}(t, \cdot, \cdot) \big)\big)_{t\geq t_0}\right]
  \La \CL \left[(u^{k,\ell}(t), U^{k,\ell}(t, \cdot,\cdot)))_{t \geq t_0}\right],
  \ee
  in law on pathspace, where the r.h.s. is
  supported on the solution of the {\em nonlinear} system
  (\ref{grocoll2}) and (\ref{ba10g}) corresponding to the initial state
  $(u(t_0), U(t_0))$. (Note that the mechanism of the limit dynamics
  does not depend on $k$ or $\ell$, but the state at $t_0$ will.)
  $\qquad \square$
  \end{proposition}

We obtain the limiting system $(u,U)$ as follows.  We specify the pair
\be{ba10d1}
(u,U) = (u(t), U(t))_{t \in \R}, \mbox{ with } (u(t), U(t)) \in \R^+\times \CM_1(\R^+\times \N),
\ee
by the (coupled) system of nonlinear {\em forward} evolution equations:

\be{grocoll2}
\frac{du(t)}{dt}= \alpha(t)(1-u(t)) u(t) -\gamma(t) u^2(t),
\ee

\bea{ba10g}
\frac{\partial U(t,dv,j)}{\partial t}&&\\
  =&&- \frac{\partial U(t,dv,j)}{\partial v} \nonumber
\\&&
+s(j-1)1_{j\ne1}U(t,dv,j-1)-sjU(t,dv,j)\nonumber
\\&&
  +\frac{d}{2}(j+1)jU(t,dv,j+1)-\frac{d}{2}j(j-1))1_{j\ne1}U(t,dv,j)\nonumber
\\&&
  +c(j+1)U(t,dv,j+1)-cjU(t,dv,j)1_{j\ne 1} \nonumber
\\&&
  -cu(t)U(t,dv,1)1_{j=1}\nonumber
\\&&
  + u(t)(\alpha(t)+\gamma(t))[1_{j\ne1}U(t,dv,j-1)-U(t,dv,j])   \nonumber
\\&&
  + \left ((1- u(t))\alpha(t)\right )1_{j=1}\cdot \delta_0(dv) \nonumber
\\&&
 -\Big(\alpha(t)(1-u(t)) -\gamma(t)u(t)\Big)\cdot U(t,dv,j).
\nonumber\eea

We have to constrain the state space to guarantee the r.h.s. above
is well-defined. Set therefore $\nu$ to be the measure on $\N$ given by
\be{agre65a}
\nu(j)=1+j^2
\ee
and consider $\R \otimes L^1 (\nu,\N)$ as a basic space for the analysis.
Then the equations (\ref{grocoll2}), (\ref{ba10g}) have the following properties.

\beP{P.UCE} {(Uniqueness of the pair $(u,U)$)}

(a) The pair of equations (\ref{grocoll2}) - (\ref{ba10g}),  given
an initial state from $\R^+ \times \CM_1(\R^+ \times \N)$
satisfying
\be{agr4}
(u(0), U(0, \cdot, \cdot)) \in \R \otimes L^1 (\N,\nu)
\ee
at time $t_0$, has  a unique solution $(u(t),U(t))_{t\geq t_0}$
with values in $\R^+ \otimes L^1_+ (\N,\nu)$,
which satisfies
$u \geq 0$ and $U(t,\R^+ \times \N)\equiv 1$.

(b) There exists a solution $(u,U)$ with time parameter $t\in
\mathbb{R}$ for every $A\in (0,\infty)$ with values in
$\R \otimes L^1(\N,\nu)$, such that
\be{S081}
u(t)e^{-\alpha t}\to A \text{  as  } t\to -\infty,\ee

\be{S082} U(t)\ttooo {\CU}(\infty).
\ee
Here $\CU(\infty)$ is the
stable age and size distribution of the CMJ-process corresponding to
the particle process
$(K_t, \zeta_t)_{t \geq 0}$
given by the McKean-Vlasov dual process $\eta$, as
defined in (\ref{ad14}).
\sm

 (c) Given any solution $(u,U)$ of
equations (\ref{grocoll2}) - (\ref{ba10g}) for  $t\in \mathbb{R}$
with values in the space $\R \otimes L^1(\N,\nu)$ satisfying \be{limsup2}
u(t) \geq 0, \quad \limsup_{t\to -\infty}e^{-\alpha t}u(t)<\infty,
\ee
we then have for $u$
\be{limsup2b}
A=\lim_{t\to -\infty}e^{-\alpha t}u(t)\ee exists and the solution satisfying
this for given $A$ is unique.  Furthermore $U$ satisfies
\be{limsup2a}
U(t)\La \CU(\infty) \mbox{ as } t \to -\infty. \qquad \square
\ee
\end{proposition}

\begin{remark}
Potential limits arising from $(u^N, U^N)$ do satisfy equation
(\ref{limsup2}).
\end{remark}

\begin{remark} Note that looking at the form of the equation we see that
a solution indexed by $\R$ remains a solution if we make a time shift.
This corresponds to the different possible values for the growth constant
$A$ in (\ref{S081}). In particular the entrance law from 0 at time
$-\infty$ is unique up to the time shift.
\end{remark}

We can now identify the behaviour of $(u^N, U^N)$ for $N \to \infty$
in terms of the early growth behaviour of the droplet of colonized sites as follows.

\beP{P.conv}{(Identification of colonization-equilibrium dynamics)}

The limits $(u^{k,\ell}(t), U^{k,\ell}(t, \cdot,\cdot)))_{t \geq
t_0},$ in (\ref{agrev60}) can be represented as the unique solution
of the nonlinear system (\ref{grocoll2})
and (\ref{ba10g}) satisfying
\be{agrev63}
\lim_{t\to -\infty}e^{-\alpha t}u^{k,\ell}(t)=W^{k,\ell},\; \lim_{t\to -\infty} U^{k,\ell}(t)
= \mathcal{U}(\infty),
\ee
with $W^{k,\ell}$ having the law of the random variable appearing as
the scaling (by $e^{-\alpha t}$) limit $t \to \infty$
of the CMJ-process $K^{k,\ell}_t$ started with $k$ particles at each of $\ell$ sites.
$\qquad \square$
\end{proposition}

We can now ask, what happens if we consider times
$\alpha^{-1} \log N + t_N$ with $t_N \to \infty$.
Then we reach a global stable state. Let
$\wt \Pi^{\iota^\ast}_{c,d,s} \Pi^{\iota^\ast}_{c,d,s}$
conditioned to be strictly positive.

\beP{P.eqpop}{(Equilibrium population)}

We have
\be{ep1}
N^{-1} u^N(\alpha^{-1} \log N + t_N) \Ntoo \iota^\ast \quad , \quad
U^N (\alpha^{-1} \log N + t_N) \Nto \wt \Pi^{\iota^\ast}_{c,d,s},
\ee
where $\iota^\ast$ satisfies the selfconsistency relation:
\be{scr1}
\suml^\infty_{k=1} k \pi^{(\iota)}_{c,d,s} (\{k\}) = \iota.
$\qad$
\ee
\end{proposition}

\section{The duality relation}
\label{s.dualrel}

In this section we relate the processes from section \ref{s.FWraremut} and \ref{s.logisticbrw}
with each other by duality and we discuss the extension of the duality to more than two types.

\subsection{A classical duality formula}
\label{ss.basicform}

The key tool in relating the two processes we have introduced is duality.
Recall the classical relation between a single Fisher-Wright  diffusion
and the Kingman coalescent.

Consider $(X_t)_{t \geq 0}$ solving
\be{ad20}
d X_t = \sqrt{d \cdot X_t (1-X_t)} dW_t
\ee
and $(D_t)_{t \geq 0}$ being the $\N$-valued death process
\be{ad21}
n \to n-1 \mbox{ at rate } d \cdot {n \choose 2}.
\ee
Then a generator calculation shows (we explain more on the background
below):
\be{ad22}
E_{X_0} [(X_t)^k] = E_k [X^{D_t}_0].
\ee

The analogous relation can be formulated for our mean-field spatial
model including mutation and selection. Define for the process
starting with $k$-particles at each of $\ell$ sites:
\be{ad28}
\Pi^{N,k,\ell}_t = u^N(t)\sum_{j=1}^\infty j U^N(t,\R^+,j),
\ee
where $(u^N(t), U^N(t,\R^+,\N))$ are given by (\ref{ba10f4}).
Then we get by a generator calculation the formula:
\be{ad23}
E [x^N_1(i,t)] = E[exp (-\frac{m}{N} \intl^t_0 \Pi^{N,1,1}_u du )],
\quad i \in \{1, \cdots, N\}.
\ee

Hence if we start initially with $x^N_1 (i,0) = 1; i=1, \cdots,N$
we see that in order to observe a mean which is less than one we need an occupation
measure of the population from the branching coalescing random walk
which is of order $N$. Due to the exponential growth this amounts
to $\Pi^{N,1,1}$ to grow up to order $N$. Since we proved the latter
behaves at times $\alpha^{-1} \log N+t$ for very negative $t$
approximately like
\be{ad24}
W exp (\alpha(\alpha^{-1} \log N+t)) = e^{\alpha t} \cdot WN,
\ee
the advantageous type emerges since
$E[x^N_1 (1, \alpha^{-1} \log N+t)]$ is as $N \to \infty$ for $t \in \R$ strictly between
0 and 1. By the analysis of moments we actually can prove the results
on the Fisher-Wright diffusion model stated in section \ref{s.FWraremut}
from the results of the coalescing particles model stated in section \ref{s.logisticbrw}.
In section 7 of \cite{DGsel} methods are developed to turn this idea into rigorous
mathematics.

\subsection{The genealogy and duality}
\label{ss.gendual}

The formula (\ref{ad22}) can be understood on a deeper level, since
the dual process can be interpreted in terms of the genealogy of the
Fisher-Wright diffusion model. This will also exhibit the role of
selection a bit better.

For that purpose the Kingman coalescent has to be viewed as a
partition-valued process. This means its states are the partitions
of the set $\{1, \cdots, k\}$ starting in $\{\{1\}, \cdots, \{k\}\}$
where partition elements coalesce at rate $d$ independent of each
other (and in the spatial case perform continuous time random walks
at rate $c$). If we consider a Fisher-Wright diffusion it can be
viewed as the limit of the Moran model taking the population size
to infinity. Pick then from the population at time $t$ exactly $k$ individuals
and look at their genealogy. This genealogy has the same law as the
genealogy of the Kingman coalescent with the genealogical distance
between two individuals defined to be the first time they are in
the same partition element. Given this genealogy we can calculate
the probability, that all individuals have the same type, one, which is
$x^k$ if $x$ is the current (time $t$) frequency of type one, in terms of the
Kingman coalescent getting the duality relation (\ref{ad22}).
See \cite{GPWmetric}, \cite{GPWmp} for more details on the
genealogical processes.

If we include selection note that if we follow the tagged sample from the
population back one of the individuals might interact with the
rest of the population (note now individuals under selection and mutation are {\em not}
anymore exchangeable) by a selection event. What results from this event
however now depends on the current types of the two involved individuals.
Therefore we have to also follow the other individual further back.
This means we expect that each time a selection event occurs our
tagged sample has to be enriched by a further particle. This is
reflected in the birth of new individuals in the dual particle model.

In order to handle multiple types and mutation the point now is that the dual process
has to be complemented by a function-valued part. The basic idea behind
this is explained in the next subsection.

\subsection{The dual for general type space}
\label{ss.dualgts}

The picture from above for two types is more subtle, if we consider
a multitype situation with more than 2 types. The basic new ideas
for this purpose in \cite{DGsel} we sketch next.

\subsubsection{A multitype model}
\label{sss.mtm}

The process considered in section \ref{s.FWraremut} can be defined
for general type space $\I$ (a subset of $[0,1]$) as
$(\CP(\I))^S$-valued process ($S$ = a finite or countable geographic space) and is called the
interacting Fleming-Viot process with selection and mutation, see \cite{DG99},
(which becomes a multitype Fisher-Wright diffusion in the case
 of finitely many types and the model of section \ref{s.FWraremut} for two types).
Then types have a fitness given by a function and mutation occurs via a
jump kernel which are denoted respectively by
\be{dg1}
\chi :\I \la \R^+, \; 0 \leq \chi \leq 1, \; M(x,dy) \mbox{ a } \I \times \I
\mbox{-probability transition kernel}.
\ee
The generator of the nonspatial Fleming-Viot process acts on monomials of order
$n$ with test function $f \in C_b(\I)$
\be{dg2}
F(x) = \int f(u_1, \cdots, u_n) x (du_1) \cdots x (du_n),
\quad x \in \CP(\I),
\ee
as follows, with setting $Q_x (du,dv) = x(du) \delta_u(dv) - x(u) x(dv)$:

\bea{dg10}
(GF)(x) &=&  s \intl_{\I}\left\{  \frac{\partial F(x)}
    {\partial x }
    (u)\left(\chi(u)-\int_{\I} \chi(w)x (dw) \right)\right\}x(du) \\
&&+\; m\intl_{\I}  \left \{ \intl_{\I}
    \frac{\partial F(x)}{\partial x } (v)M(u,dv)
   -  \frac{\partial F(x)}{\partial x }(u)\right \}
      x (du) \nonumber \\
&&+ \;d\intl_{\I} \intl_{\I}
     \frac{\partial ^2F(x)}{\partial x \partial x }
     (u,v)Q_{x }(du,dv)  \Bigg ], \qquad x \in
	 \CP(\I). \nonumber
\eea
In the spatial case a corresponding drift term from migration appears as well:
\be{dg10b}
\suml_{i,j\in S} a(i,j) \int \big(\frac{\partial F(x)}{\partial x_j} (u)
- \frac{\partial F(x)}{\partial x_i}(u)\big) x_i(du).
\ee

In \cite{DK} a dual for a Fleming-Viot process
 with mutation and selection is introduced in order to show that the process is well-defined by its martingale problem.  In \cite{DGsel} a new dual was developed which makes possible the study of the long-time behavior and the genealogy of the system.
This we explain now.
 In order to introduce the main ideas we first consider the case with no mutation and the
 special case $N=1$ and then the general case.

\subsubsection{The dual with selection}
\label{sss.dualsel}

Let now the mutation rate be zero and let $|S|=1$, so that the state space
is $\CP(\I)$ and we only have resampling and selection.
Consider the class of functions
\be{3.2}
 F((n, f), x) : = \left [
 \intl_{\I} \dots \intl_{\I} f(u_1,\dots,u_n )
 x(du_1) \dots x(du_n)\right ], \ee
for all $n \in \N$,  $f \in L_\infty((\I)^n, \R)$ and $x\in\mathcal{P}(\mathbb{I})$.
Note that given a random probability measure
$X$ on $\mathcal{P}(\mathbb{I})$, the collection
\be{dg3}
\{E[F((n,f),X)], n\in\N,  f \in L_\infty((\I)^n, \R)\}
\ee
uniquely characterizes the probability law of $X$.

The class given above contains the functions which we will use
for our dual representations.
In fact it suffices to take smaller, more convenient sets of test functions $f$:
\be{productform}
  f(u_1,\dots,u_n )=\prod_{i=1}^n f_{i}(u_i),
  \quad f_i \in L_\infty (\I). \ee
If we consider the case when $\mathbb{I}$ is
finite, it would even suffice to take functions
  \be{3.2b}
  f_{j}(u)=1_{j}(u) \mbox{ or } f_j(u)=1_{A_j}(u),
\ee
where $1_{j}$ is the indicator function of $j\in\mathbb{I}$. \bi

The \emph{function-valued dual processes}
$(\eta_t, \CF_t)_{t \geq 0}$ and
$(\eta_t, \CF^+_t)_{t \geq 0}$ are
constructed from the following four ingredients:

\begin{itemize}
\item $N_t$ the number of {\em individuals} present in the dual process which
is a non-decreasing $\N$-valued process with $N_0=n$, the number of initially tagged individuals,
and $N_t\geq n$,
\item $\zeta_t =\{1,\cdots,N_t\}$ is an {\em ordered particle system},
\item   $\pi_t$: is a  {\em partition} $(\pi^1_t,\cdots,\pi^{|\pi_t|}_t)$ of $
\zeta_t$, i.e. an ordered family of subsets, where the
{\em index} of a partition element is the smallest element of the partition element,\\

\item $\CF_t$ is for given states of $\pi_t, \zeta_t,N_t$ a function in
$L_\infty (\I^{|\pi_t|})$ which is obtained from a function in
$L_\infty(\I^{N_t})$ by setting variables equal which are corresponding to  one
and the same partition element {and $\CF_t$ changes further driven
selection (see below)}.
\end{itemize}

\beD{D.Fuvadure} {(Evolution of $(\eta, \CF)$ and $(\eta, \CF^+)$)}

(a) The process $\eta$ is driven by coalescence at rate $d$ of every pair
of partition elements and by the birth of a new individual at rate $s$,
which forms its own partition element.

(b) Conditioned on the process $\eta$ the evolution of $\CF$ is as follows.

\begin{itemize}
\item The coalescence mechanism:
 If a coalescence of two partition elements occurs,
 then the corresponding variables of $\CF_t$ are set equal to the
 variable indexing the partition element,
 i.e. for $\CF_{t-}=g$ we have the transition
 (here $\wh u_j$ denotes an omitted variable)
\bea{ga3}
 g(u_1,\cdots,u_i,\cdots,u_j,\cdots,u_m)
 & \longrightarrow & \wh g(u_1,\cdots,\wh u_j,\cdots,u_{m}) \\
 && =  g(u_1,\cdots,u_i,\cdots,u_{j-1}, u_i, u_{j+1},\cdots,u_{m}), \nonumber
\eea
so that the function changes from an element of $L_\infty(\I^m)$ to one of
$L_\infty(\I^{m-1}).$
\item The selection mechanisms:
\begin{itemize}
\item \underline{Feynman-Kac.}  For $\CF_t$, if a birth occurs in the process $(\zeta_s)$ due to the
partition element  to which the element $i$ of the basic set belongs, then
for $\CF_{t-}=g$ the following transition occurs from an element in
$L_\infty((\I)^m)$ to elements in $L_\infty((\I)^{m+1})$:
\be{3.6}
g(u_1,\dots,u_m) \longrightarrow \chi(u_i)g(u_1,\dots,u_m) -
\chi(u_{m+1}) g(u_1,\dots,u_m). \ee
\item \underline{Non-negative.} For $\CF^+_t$ the transition (\ref{3.6}) is replaced by (provided that
$\chi$ satisfies $0 \leq \chi \leq 1$):
\be{3.6b}
g(u_1,\cdots,u_m) \longrightarrow \wh g(u_1,\dots,u_{m+1})=(\chi
(u_i) + (1-\chi (u_{m+1}))) g(u_1,\cdots,u_m). $\qad$
\ee
\end{itemize}
\end{itemize}
\end{definition}

Now we can obtain two different duality relations, the first below working
in all cases and a second below working in a large subclass of models.

\beP{P3.1} \quad (Duality relation - signed  with Feynman-Kac dual)\\
Let $(X_t)_{t \ge 0}$ be a solution of the Fisher-Wright martingale
with finite type space $\mathbb{I}$, fitness function $\chi$
and selection rate $s$ with $X_0=x\in\mathcal{P}(\mathbb{I})$. Let $\CF_0=f \in L_\infty((\I)^n)$ for some $n \in \N$.

Assume that $t_0$ is such that:
\be{3.1.1}
\label{expcondition} E\left(\exp\left(s\int_0^{t_0}|\pi_r|dr\right
)\right)<\infty.\end{equation}

Then for $0\leq t\leq t_0$, $(\eta_t, \CF_t)$ is the Feynman-Kac
dual of $(X_t)$, that is:
\bea{3.8}
\lefteqn{E[F(( \eta_0,f),
X_t)] =  E_{( \eta_0, \CF_0)}\Bigg \{ \Bigg [ \exp(s \intl^t_0
|\pi_r|dr) \Bigg ] \cdot  }\\
[1ex] & & \hspace*{3cm}  \cdot \Bigg [\intl_{\I} \dots \intl_{\I}
\CF_t(u_1,\dots,u_{|\pi_t|}) x (du_1)\dots
x(du_{|\pi_t|}) \Bigg ] \Bigg\}, \nonumber%\\
% & &  = E_{( \pi_0, \CF_0)} \cdot \Bigg [\intl_{\I} \dots \intl_{\I}
%\CF^+_t(u_1,\dots,u_{|\pi_t|}) x (du_1)\dots
%x(du_{|\pi_t|}) \Bigg ] \Bigg\}, \nonumber
 \eea
where the initial state $(\eta_0, \CF_0)$ is for $n \in \N$ chosen given by

\bea{3.9}
 \pi_0 & =&  [(\{1\},\{2\},\dots,\{n\})],
 \\
 \CF_0 & =&  f \in C_b(\I^n). \qquad \square \nonumber \eea
 \end{proposition}

 \begin{remark} In \cite{DG99} it was shown that there exists $t_0 > 0$ for
 which (\ref{3.1.1}) is satisfied.
 \end{remark}
 \bi

The disadvantage of the dual above is the exponential term together with the signed function.  This involves the interplay of a cancelation effect and the exponential growth factor which
is often hard to analyse as $t \to \infty$.  The key observation is that if the fitness function $\chi$
is a bounded function
 then we can obtain the following duality relation that
 does {\em not} involve a Feynman-Kac factor and preserves the positivity of functions:

\beP{P.orddual}(Duality relation - non-negative)

With the notation and assumptions as in Proposition \ref{P3.1}
(except with (\ref{3.6}) replaced by (\ref{3.6b})) we get for $\chi$ with
\be{3.9.0}
0 \leq \chi \leq 1
\ee
that for all $t \in [0,\infty)$ we have:
\be{3.9.1}
E[F (( \eta_0,f),X_t)]=E_{(\eta_0,\CF_0^+)}
\left [ \intl_\I \cdots
\intl_\I  \CF^+_t (u_1, \cdots,u_{|\pi_t|}) x_{1}
(du_1) \cdots x_{ (|\pi_t|)} (du_{|\pi_t|}) \right ]. \ee
 Moreover, $\CF^+_t$ is always non-negative if $\CF^+_0$ is.\hfill $\square$
 \end{proposition}

\subsubsection{The dual with migration, selection and mutation}
\label{sss.dualsmraremut}

We now consider the case of (\ref{angr2}), (\ref{angr3}) including selection and migration but for the moment setting $m=0$.
The partition elements of the dual now have {\em locations} given by
\be{dg4}
\xi_t: \pi_t \longrightarrow \Omega^{|\pi_t|}_N.
\ee

The corresponding additional {\em dual migration} dynamics is as follows.
At rate $c$ each particle can jump from its current location to a randomly chosen point
in $\{1,\dots,N\}$.  In the limiting case ($N \to \infty$)
the particle always migrates to an empty site which for convenience we can take to be the smallest unoccupied site $n\in\N$.
This then results in precisely the logistic branching particle model respectively
the Crump-Mode-Jagers process described in section \ref{s.logisticbrw}.

The duality relation is now given by
\be{3.9.2}
E[F (( \eta_0,f),X_t)]=E_{(\eta_0,\CF_0^+)}
\left [ \intl_\I \cdots
\intl_\I  \CF^+_t (u_1, \cdots,u_{|\pi_t|}) x_{\xi_t(1)}
(du_1) \cdots x_{\xi_t (|\pi_t|)} (du_{|\pi_t|}) \right ],
\ee
where $X_0=(x_1,\dots,x_n),\;x_i\in\mathcal{P}(\mathbb{I}),\;
i=1,\dots,N$, $\eta_0=(\zeta_0,\pi_0,\xi_0)$.

The mutation can now be incorporated by adding a corresponding transition of
$\CF$, at rate $m$ for every variable
the operation acting on $g \in (L^\infty (\I))^n$ by
\be{dg5}
g(u_1, u_2, \cdots, u_n) \la \int g (u_1, \cdots, v, u_{i+1}, \dots, u_n)
M(u_i, dv).
\ee
Then the relation (\ref{3.9.1}) holds with mutation.

%If we now have to handle also the mutation and the selection the
%picture becomes more complicated. However in \cite{DGsel} a
%duality relation is developed which allows to handle multitype
%situations with mutation and selection in such a way that
%a duality relation between the Fisher-Wright model and a
%dual process can be established where the particle model
%we introduced in section \ref{s.logisticbrw} is one key ingredient.
%The first question which arises is why does the dual have births besides
%the coalescence.

{\bf Example of application}
Now consider the case with $m>0$, $\mathbb{I}=\{1,2\}$ as in section 1 and
$\CF_0 =1_1$.
This effect of rare mutation from type 1 to type 2 results in the dual in the transition
\be{dg6}
1_1\to 0\quad \text{at rate  }\frac{m}{N}.\ee
Returning to our duality relation we now see that the particles in our
process $\eta$ stand for possible individuals, represented by the factor $1_1$ which could by mutation
represent a possible line through which the advantageous type can enter
the population.  More precisely, if $\mathcal{F}_0=1_1$, then we can check that as time increases
we have an increasing number, $\Pi^{N,1,1}_u$ of such factors and any of these can undergo the rare mutation transition.

Therefore
\be{dg6b}
V^N_t = \frac{m}{N} \intl^t_0 \Pi^{N,1,1}_u du
\ee
represents the {\em hazard function} for a rare mutation to occur and therefore
$(1-\exp (-V_t^N))$ the mean of $x^N_2(t)$. Hence the growth of
the spatial intensity in the logistic branching particle model
relates to the emergence of the rare mutant population.

As a simple application of this dual we can now verify (\ref{ad2}), (\ref{ad3}) when there is only one site.  In the case $d=0$ we have a linear birth process with birth rate $s$ and therefore $V^N_t$  grows like $e^{st}$ and therefore the first rare mutation occurs after a time of order $O(\frac{\log L}{s})$.  On the other hand when $d>0$ the process does not grow indefinitely but approaches an equilibrium.  In this case $V^N_t$ grows only in a linear fashion and therefore requires time of order $O(L)$.

We note  that since in the process, respectively the dual, time has
to be read forward, respectively backward, from $t$, the rare mutation jumps
occuring in times $\alpha^{-1} \log N+t$ for the dual after the
expansion of the population up to time $\alpha^{-1} \log N$,
correspond to rare mutations in the original process
occurring at times between 0 and $t$ and then
growing till the takeover after time $\alpha^{-1} \log N$.

\subsection{Outlook on set-valued duals}
\label{ss.outlooksvd}

If we consider $| \I | <\infty$ and use for $\CF^+_0$ {\em products of
indicator functions} we obtain under the evolution {\em sums of products
of indicator functions} where the dynamic of the different summands
is coupled by the transition occuring in the underlying process $\eta$.
Here we briefly  describe the main idea how to describe the dual
based on a set-valued process but where we introduce the order of
individuals and we use a change in the coupling between
summands.

To explain the main idea we again consider the case $\mathbb{I}=\{1,2\}$, $N=1$.
Now consider the case
\be{dg7}
\mathcal{F}_0= 1_2.\ee
Then at the random time $\tau$ at which one selection operation occurs we have

\be{dg8}
\mathcal{F}_\tau= 1 \otimes 1_2 +1_1\otimes 1_2.\ee
We can now regard this as defining a subset of $\mathbb{I}^2$.  If we now couple the transitions in the two summands
differently, namely we write instead of (\ref{dg8})
(this we can do since the dual expression depends
only on the marginal law of the summands not the joint law)
\be{dg9}
1_2 \otimes 1 +1_1\otimes 1_2,
\ee
then
we can ensure that the summands continue to correspond to {\em disjoint}
subsets of $\mathbb{I}^\N$ and therefore we obtain a dual process with values in subsets
of $\mathbb{I}^\N$.

A much more complicated set-valued dual can be constructed for the general multitype Fisher-Wright diffusion with mutation, selection and migration.
The key point is to first take the
non-negative function-valued dual driven by the particle system we introduced and then to introduce the
{\em order of factors} and an appropriate {\em coupling} of its decomposition into a set of summands.
This allows us to obtain a duality relation for general finite type space and
additive selection with a bounded fitness function.

This {duality relation} can be  interpreted
in terms of the {\em ancestral lines of a tagged sample} of individuals picked from the time-$t$
population, since the dual now gives a decomposition in disjoint events for the
ancestral lines and genealogical tree for a tagged sample of $n$ individuals
from the time $t$ population from which we can read off current types and
genealogical distance of the tagged sample.


\begin{thebibliography}{9999999}

\bibitem[DGsel]{DGsel}
D. Dawson and A. Greven:
{\em On the effects of migration in spatial Fleming-Viot models with
selection and mutation}, in preparation 2010.

\bibitem[Bu]{Bu01} R. B\"urger (2001).  The Mathematical Theory of Selection,
Recombination, and Mutation, Wiley.

\bibitem[DG99]{DG99} D. A. Dawson and A. Greven (1999),
\emph{Hierarchically interacting Fleming-Viot processes with
selection and mutation: Multiple space time scale analysis and
quasi equilibria}, Electronic Journal of Probability, Vol. {\bf
4}, paper no. 4, pages 1-81.


\bibitem[D]{D} D. A. Dawson (1993), {\it Measure-valued Markov
Processes.} In: \'Ecole d'\'Et\'e de Probabilit\'es de Saint Flour
XXI, Lecture Notes in Mathematics {\bf 1541}, 1-261,
Springer-Verlag.

\bibitem[DG2]{DG2} D. A. Dawson and A. Greven (1993b),
{\it Multiple time scale analysis of interacting diffusions.}
Probab. Theory Rel. Fields {\bf 95}, 467-508.

\bibitem[DG3]{DG3} D. A. Dawson and A. Greven (1993c),
{\it Hierarchical models of interacting diffusions: multiple time
scale phenomena. Phase transition and pattern of
cluster-formation.} Probab. Theory Rel. Fields, vol. {\bf 96},
435-473.

\bibitem[DK]{DK} D.A. Dawson and T.G. Kurtz (1982). Applications of duality to measure-valued diffusions,
Springer Lecture Notes in Control and Inf. Sci. 42, 177-191.

\bibitem[EK4]{EK4}
S. N. Ethier and T. G. Kurtz (1994), {\it Convergence to
Fleming-Viot processes in the weak atomic topology}, Stochastic
Process, Appl. {\bf 54}, 1-27.



\bibitem[Gar]{Gar} J. G{\"a}rtner (1988), {\it On the McKean-Vlasov limit for interacting diffusions}.
Math. Nachr. {\bf 137}, 197-248.

\bibitem[GPWmp]{GPWmp}
   A. Greven,  P. Pfaffelhuber and A. Winter,
	 {\em Tree-valued resampling dynamics: Martingale Problems and applications},
	      submitted to PTRF 2010.

\bibitem[GPWmetric]{GPWmetric}
   A. Greven,  P. Pfaffelhuber and A. Winter (2009),
   {\em Convergence in distribution of random metric measure
    spaces (The $\Lambda$-coalescent measure tree)},
 PTRF, Vol. 145, issue 1, 285 ff.


\bibitem[Hu]{Hu} M. Hutzenthaler (2009), \emph{The virgin island model}, Electr. J. Probab. 14, 1117-1161.

\bibitem[J92]{J92} P. Jagers (1992), {\it Stability and
instability in population dynamics}, J. Appl. Probab. 29, 770-780.



\bibitem[JN]{JN} P. Jagers and O. Nerman (1984), {\em The growth and
composition of branching populations}, Adv. in Appl. Probab. 16,
221-259.

\bibitem [Mc]{Mc} H. P. McKean, Jr. (1966), A class of Markov processes associated with nonlinear parabolic
equations, Proc. N.A .S., U.S.A. 56:1907-19

\bibitem[N]{N} O. Nerman (1981), {\em On the convergence of
supercritical general (C-M-J) branching processes}, Zeitschrift f. Wahrscheinlichkeitsth. verw.
Gebiete, 57, 365-395.




\bibitem[PY]{PY} J. Pitman and M. Yor (1982), {\em A decomposition of Bessel
bridges}, Z.  Wahr. verw. Geb. 59, 425-457.

\bibitem[RW]{RW} L.C.G. Rogers and D. Williams (1987). Diffusions,
Markov processes and martingales, Vol. 2, Wiley.


\bibitem[Schirm10]{Schirm10} F. Schirmeier (2010),
{\em A spatial population model in separating time windows},
master thesis, Department Mathematik, Erlangen.


\end{thebibliography}
\end{document}